\documentclass [11pt]{amsart}
\usepackage {amsmath, amssymb, amscd, mathrsfs, hyperref, slashed, enumerate, url}
\usepackage[text={6in,8.5in},centering,letterpaper,dvips]{geometry}
\usepackage[all,cmtip]{xy}

\theoremstyle{remark}

\DeclareFontFamily{U}{mathx}{\hyphenchar\font45}
\DeclareFontShape{U}{mathx}{m}{n}{
      <5> <6> <7> <8> <9> <10>
      <10.95> <12> <14.4> <17.28> <20.74> <24.88>
      mathx10
      }{}
\DeclareSymbolFont{mathx}{U}{mathx}{m}{n}
\DeclareFontSubstitution{U}{mathx}{m}{n}
\DeclareMathAccent{\widecheck}{0}{mathx}{"71}

\def\polhk#1{\setbox0=\hbox{#1}{\ooalign{\hidewidth
    \lower1.5ex\hbox{`}\hidewidth\crcr\unhbox0}}}

\def\del {\partial}

\begin{document}

\title{Errata to the article ``Seiberg-Witten-Floer stable homotopy of three-manifolds with $b_1=0$''}
\author[Ciprian Manolescu]{Ciprian Manolescu}
\address {Department of Mathematics, UCLA, 520 Portola Plaza\\ 
Los Angeles, CA 90095}
\email {cm@math.ucla.edu}
\maketitle

There is a significant error at the top of p.921 in the paper. The sentence ``Changing everything by a gauge, we can assume without loss of generality that $i^*(\hat a) \in \operatorname{ker } d^*$'' is not correct. There is no gauge freedom available on $X$, because on p.916 we already fixed the gauge by imposing the Coulomb-Neumann condition on forms, $\hat{a} \in \Omega^1_g(X)$. 

Furthermore, if we do not assume that $i^*(\hat a) \in \operatorname{ker } d^*$, then the estimates on the second term on the right hand side of Equation (17) do not go through. Precisely, in the paragraph on p.921 starting with ``Similarly one can show that \dots,'' instead of $db_n \to 0$ in $L^2_k$ and $b_n \to 0$ in $L^2_{k+1}$ we would have $db_n \to db$ in $L^2_k$ and $b_n \to b$ in $L^2_{k+1}$, where $i^*(\hat x) = (a+db, \phi)$ and $a \in \operatorname{ker } d^*$. Knowing that $b_n \to b$ and $p^{\mu_n}_{\lambda_n} (a_n, e^{ib_n} \phi_n) \to (a, e^{ib}\phi)$, we would like to deduce that $p^0(a_n, \phi_n) \to p^0(a, \phi)$. (Here, all limits are in $L^2_{k+1}$.) By hypothesis, we also know that the $L^2_{k+1}$ norms of $(a_n, \phi_n)$ are bounded, and that $(a_n, \phi_n) \in V_{\lambda_n}$. Thus, $p^0(a_n, \phi_n) = p^0_{\lambda_n} (a_n, \phi_n)$. We have
\begin{align*}
  \|p^0 (a_n, \phi_n) - p^0 (a, \phi) \| &= \|  p^0_{\lambda_n} (a_n, \phi_n) - p^0(a, \phi) \| \\
  & \leq \|  p^{0}_{\lambda_n} (a_n, e^{ib_n-ib} \phi_n) - p^0(a, \phi) \| + \| p^0_{\lambda_n}\bigl(0, (e^{ib_n - ib} - 1) \phi_n\bigr)\|,
  \end{align*}
where all norms are $L^2_{k+1}$. Since $b_n \to b$ and $\|\phi_n\|_{L^2_{k+1}}$ is bounded, using the Sobolev multiplication $L^2_{k+1} \times L^2_{k+1} \to L^2_{k+1}$ we get that the second term in the last expression above converges to $0$. 
 If multiplication by $e^{ib}$ commuted with the projection $p^{\mu_n}_{\lambda_n}$, from $p^{\mu_n}_{\lambda_n} (a_n, e^{ib_n} \phi_n) \to (a, e^{ib}\phi)$ we would get that $p^{\mu_n}_{\lambda_n} (a_n, e^{ib_n-ib} \phi_n) \to (a, \phi)$, and then (applying $p^0$) the first term would converge as well. It would then follow that $p^0(a_n, \phi_n) \to p^0(a, \phi)$, as desired.
 
This argument works for $b=0$ but fails in general, because multiplication by $e^{ib}$ does not commute with $p^{\mu_n}_{\lambda_n}$. The origin of the problem is that the nonlinear map $C^{\mu}$ defined on p.917 is not compact.

The simplest way to fix this issue is to replace the Coulomb-Neumann condition by a double Coulomb condition. This approach is the subject of Khandhawit's paper \cite{Khandhawit}. We sketch the argument here, and refer to \cite{Khandhawit} for more details. 

On p.916, when we define $\Omega^1_g(X)$, instead of the condition $\hat{a}|_{\del X} (\nu) =0$ we impose a boundary Coulomb condition, $i^*(\hat a) \in \ker(d^*)$. We also ask that the integral of $\hat{a}|_{Y_i} (\nu)$ is zero on each connected component $Y_i \subseteq \del X$. (This is automatic when $\del X$ is connected.) The new gauge condition satisfies a Fredholm property similar to Proposition 5; see \cite[Proposition 2]{Khandhawit}. Moreover, the nonlinear map $C^\mu$ from p.917 is now compact, and we can delete $\operatorname{pr}_{\operatorname{ker } d^*}$ from the second term on the right hand side of Equation (17) on p.921. Then, it is easy to show that this term converges to zero. A new difficulty appears in the argument at the top of p.922, when we glue a half-trajectory on $[0, \infty) \times Y$ with a monopole on $X$ that may have a non-trivial $dt$ component on the boundary. Nevertheless, the gluing can be done after changing the half-trajectory on $[0, \infty) \times Y$ by a suitable gauge transformation; see \cite[Corollary 2]{Khandhawit}. 

\medskip
There were a few other minor errors in the article:

\medskip
\begin{enumerate}
\item On p.898, the metric $\tilde g$ on $V$ was defined by the formula
$$ \| (b, \psi) \|_{\tilde g} = \| (b, \psi) + (-id\xi, i \xi \phi)\|_{L^2},$$
measuring the norm of the projection of $(b, \psi)$ to the local Coulomb slice at $(a, \phi)$. However, this formula does not yield a non-degenerate metric. There is still a residual $S^1$ gauge action on $V$, and the vectors tangent to the $S^1$-orbits, such as $(0, i\phi)$, would have length zero. We can correct this by adding a circular projection term, given by the square of the inner product with $(0, i\phi)$. Precisely, we set:
$$ \| (b, \psi) \|^2_{\tilde g} = \| (b, \psi) + (-id\xi, i \xi \phi)\|^2_{L^2} + \bigl( \operatorname{Re} \langle i\phi, \psi \rangle \bigr) ^2.$$
Since the gradient of the $CSD$ functional is perpendicular to the $S^1$-orbits, it is still true that the trajectories of the $\tilde g$-gradient of $CSD|_V$ are the Coulomb projections of the trajectories of $CSD$ on $i\Omega^1(Y) \oplus \Gamma(W_0)$.

\medskip
\item
In the middle of p.907, when we define the desuspension of $X$ by $E$ in the category $\mathfrak{C}_0$, the alternative definition as $\Omega^E X$ is incorrect. The correct definition is the one given in the previous line, $\Sigma^{-E}X = (E^+ \wedge X, 2 \dim E, 0)$. In general, $\Sigma^{-E}X$ and $\Omega^E X$ may not even have the same homology, so they are not isomorphic in $\mathfrak{C}_0$.

\medskip
\item
At the bottom of p.917, the set  $\tilde K$ should be the preimage of $B(U_n , \epsilon_n) \times V^{\mu}_{\lambda}$ under the map $\operatorname{pr}_{U_n \times V^{\mu}_{\lambda}} \mathit{SW}^{\mu}$, not under the linear map $L^\mu$.

\medskip
\item Lemma 4 on p.918 is incorrect as stated. There can be trajectories that start outside the ball $\overline{B(2R)}$, go inside $\overline{B(2R)}$ at some time $t_0$, and converge to a point in $B(R)$. For the lemma to be true, we need an additional hypothesis, that $x_n(t_0)$ is the restriction of an approximate Seiberg-Witten solution on a compact $4$-manifold $X$ with boundary $Y$. An argument of this type is used in the proof of Lemma 2 in \cite{Khandhawit}. 

\medskip
\item At the top of p.924, the proof that  the class $\Psi$ is independent of the choices made in the construction was incomplete. One needs to show independence of the index pair $(N, L)$ chosen in Theorem 4. This is done by Khandhawit in Proposition 5 from \cite[Appendix A]{Khandhawit}.
\end{enumerate}

\medskip
{\bf Acknowledgements.} I would like to thank Tirasan Khandhawit for pointing out two of the errors  discussed above, as well as for fixing them in his paper \cite{Khandhawit}. I am also grateful to Mikio Furuta, Tye Lidman and Jianfeng Lin for discussions that led to finding the other issues.

\bibliography{biblio}
\bibliographystyle{amsplain}

\end{document}